\documentclass{amsart}
\usepackage{amssymb}
\usepackage{graphicx}
\usepackage[english]{babel}
\usepackage{subcaption}

\usepackage[pdftitle={Inverse Parabolic},
  pdfauthor={Doghonay, Engblom, Kreiss},
  pdffitwindow=true,
  breaklinks=true,
  colorlinks=true,
  urlcolor=blue,
  linkcolor=red,
  citecolor=red,
  anchorcolor=red]{hyperref}
\usepackage[numbers]{natbib}

\numberwithin{equation}{section}
\numberwithin{table}{section}
\numberwithin{figure}{section}

\newtheorem{thm}{Theorem}
\newtheorem{lem}{Lemma}
\newtheorem{rmk}{Remark}
\newtheorem{pf}{Proof}

\begin{document}

\title[Inverse heat source problems]{Efficient low rank approximations for parabolic control problems with unknown heat source}

\author[D.~Arjmand]{Doghonay Arjmand}
\author[M.~Ashyraliyev]{Maksat Ashyraliyev}
\email{doghonay.arjmand@it.uu.se,maksat.ashyraliyev@mdu.se}

\address{Division of Scientific Computing \\
  Department of Information Technology \\
  Uppsala University \\
  SE-751 05 Uppsala, Sweden.}

\address{Division of Mathematics and Physics \\\
  Mälardalen University \\
  SE-721 23 Västerås,  Sweden.}





\begin{abstract}
 An inverse problem of finding an unknown heat source for a class of linear parabolic equations is considered. Such problems can typically be converted to a direct problem with non-local conditions in time instead of an initial value problem. Standard ways of solving these non-local problems include direct temporal and spatial discretization as well as the shooting method, which may be computationally expensive in higher dimensions. In the present article, we present approaches based on low-rank approximation via Arnoldi algorithm to bypass the computational limitations of the mentioned classical methods. Regardless of the dimension of the problem, we prove that the Arnoldi approach can be effectively used to turn the inverse problem into a simple initial value problem at the cost of only computing one-dimensional matrix functions while still retaining the same accuracy as the classical approaches. Numerical results in dimensions $d=1,2,3$ are provided to validate the theoretical findings and to demonstrate the efficiency of the method for growing dimensions.
\end{abstract}

\maketitle

\section{Introduction} \label{sec:intro}

The theory of inverse problems for differential equations has been extensively developed to tackle various problems in applied sciences and engineering. Typical applications turning this field into a popular area of research include geological exploration, medical diagnostics, and predictive material science. In a direct problem, typically a differential equation is given and its solution is determined uniquely from imposed initial and/or boundary conditions. In an inverse problem, on the contrary, a differential equation itself contains unknown data/parameters. For example, a differential equation may involve unknown coefficients or source terms, which need to be determined from a set of observations or measurements of the solution. In the study of inverse problems, the lack of knowledge about the data in the model problem is then compensated by introducing some additional conditions into the problem, which may eventually result in a well-posed mathematical model. For the general theory of inverse problems for differential equations we refer the interested reader to \cite{Prilepko_book,Isakov_book,Kabanikhin_book} and the references therein.


{\color{black} Motivated by real applications, inverse problems for diffusion equations can be categorized into two: a) determining the coefficient of the equation, see e.g., \cite{Kimura1993,Chen2006,Yang2008} b) determining the heat source. Latter is closely related to the goal of the present article. In particular, we are interested in the inverse problem of finding the pair of solutions $\{ p,u\}$ such that 
\begin{equation} \label{Eqn_u}
\left\{
\begin{array}{l}
\partial_{t} u(t,x) - \nabla \cdot \left( a(x) \nabla u(t,x) \right) = f(t,x) + p(x), ~ (t,x) \in (0,T)\times K, \\ \rule{0cm}{5mm}
u(t,x) = 0, \quad (t,x) \in [0,T]\times \partial K, \\ \rule{0cm}{5mm}
u(0,x) = 0, \quad u(T,x) = \varphi(x), \quad x\in K,
\end{array}
\right. 
\end{equation}
where $K:= (0,1)^{d}$ is an open unit cube with a boundary $\partial K$ and  $a(x) \in \mathbb{R}^{d\times d}, \forall x\in K$ is a symmetric matrix function satisfying
\begin{equation}
0 < a_{\mathrm{min}} |\zeta|^{2} \leq \zeta^T a(x) \zeta \leq a_{\max} |\zeta|^{2} , \quad \forall x \in K ~ \mathrm{and} ~ \zeta \in \mathbb{R}^{d}, 
\end{equation}
and $f$ and $\varphi$ are sufficiently smooth known functions. Problem \eqref{Eqn_u} is highly relevant from an application point of view. Consider, for example, a domain $K$, for which there already exists a known (and practically difficult to remove) background heat source $f$. Suppose that our objective is to keep the temperature distribution at a specific time $T$ in a desired manner. The problem is then to find the control data $p$ accounting for yet another heat source, to be fine-tuned over $K$, such that the final temperature requirements are met.

The problem of determining the heat source in a parabolic equation has been a topic of research for the past two decades and is still a popular area of research due to the need to improve the efficiency of numerical methods linked to the proposed strategies. A class of algorithms rely on obtaining a transformed equation, where the unknown heat source is no longer present, and the resulting system is then solved with a direct numerical simulation in time and space, see e.g., \cite{Dehghan2001,Dehghan2009,Yang2011,Ashyralyev2012}. This idea typically leads to either a large linear system of equations to solve due to the urge to treat the additional temporal dimension similar to spatial dimensions rather than using a standard time-stepping method or alternatively exploiting iterative solution methods, which may be computationally expensive.  More standard ways of tackling this problem is based on re-formulation of the original problem as a minimization problem, \cite{MahmoodEtal_21,ZuiChaEtal_2022,Biccarietal_22,Huang1992,Liu2008}, which eventually needs an efficient iterative solution method to find the minimizer. There are other strategies which are limited as they assume a specific structure for the parabolic operator (e.g., constant coefficients), \cite{Cheng2020}, or they have limitations due to dimension \cite{beatriz2020solving}.  }

We may rewrite the above problem in abstract form as
\begin{equation} \label{Eqn_u_abstract}
\left\{
\begin{array}{l}
\displaystyle \frac{d}{dt}u(t) +Au(t) = f(t) + p, \quad t\in (0,T), \\ \rule{0cm}{6mm}
u(0) = 0, \quad u(T) = \varphi, 
\end{array}
\right. 
\end{equation}
where $u(t)$, $t \in [0,T]$ and $f(t)$, $t \in (0,T)$ are vector valued functions in a suitable Hilbert space $H$. Similarly, $p,\varphi \in H$ and $A: H\to H$ is the corresponding abstract operator, where the boundary conditions are incorporated in $A$. In general, problem \eqref{Eqn_u_abstract} is ill-posed. However, for sufficiently smooth data its well-posedness has been shown in the literature. The unique solvability of inverse source identification problem \eqref{Eqn_u_abstract} in an arbitrary Banach space $E$ with linear unbounded operator $A$ was established in \cite{Eidelman1991}. The stability estimates for the solution of problem \eqref{Eqn_u_abstract} were obtained in \cite{Ashyralyev2011}. {\color{black} 
Related to the theoretical foundations of the present work is also the well-posedness of an inverse parabolic problem with an unknown source term \cite{Choulli1999,Hasanov2013}; see also \cite{Cao1997,Amann2005} for regularity conditions for the controllability of final time over-determination parabolic equations with a time and space dependent right hand side, and with superlinear nonlinearities.}

{\color{black} The starting point that we follow is similar to \cite{Ashyralyev2012} which relies on  eliminating} the unknown $p$ in \eqref{Eqn_u_abstract} by introducing a new abstract function $v$ as
\[ 
v(t) = u(t) - A^{-1} p, \quad t \in [0,T]. 
\]
Indeed, it follows that $\frac{d}{dt}v(t)=\frac{d}{dt}u(t)$, $Av(t) = A u(t) - p$, $v(T) = \varphi - A^{-1} p$, and $v(0) = -A^{-1} p$. Therefore, the problem \eqref{Eqn_u_abstract} is equivalent to the following problem: 
\begin{equation} \label{Eqn_v}
\left\{
\begin{array}{l}
\displaystyle \frac{d}{dt}v(t) + Av(t) = f(t), \quad t\in(0,T), \\ \rule{0cm}{6mm} 
v(0) = v(T) - \varphi.
\end{array}
\right.
\end{equation}
Note that having the solution $v(t)$, $t\in[0,T]$ of problem \eqref{Eqn_v} allows us to obtain the solutions of initial problem \eqref{Eqn_u_abstract} directly as follows:
\[ 
p=-Av(0), \quad u(t)=v(t)-v(0), ~ t\in (0,T). 
\]
The problem \eqref{Eqn_v} is direct since it does not involve an unknown source term $p$. However, it has a non-local nature due to the dependency of $v(0)$ on $v(T)$. 

\begin{rmk}
{\color{black} Apart from the well-posedness of the problem \eqref{Eqn_u_abstract}, the stability and uniqueness of the non-local problem \eqref{Eqn_v} can be established using \cite{Ashyralyevetall2022,Starovoitov_21},  where the authors consider a more general nonlocal condition that also applies to the setting of \eqref{Eqn_v}.} 
\end{rmk}

{\color{black} The main aim of this article is to develop numerical methods, based on low rank approximations, which bypasses the limitations of the existing methods with respect to dimension and structural properties of the operator, see also \cite{MARTINVAQUERO2019166,Yanping97,Ashyralyevetall2022_b} for fully discrete approaches available in the literature. To put everything into a framework, assume a spatial discretization of the problem \eqref{Eqn_v} in the form} 
\begin{equation}
\left\{
\begin{array}{l}
\displaystyle \frac{d}{dt}v_h(t)+A_{h}v_h(t)=f_{h}(t), \quad 0<t<T, \\ \rule{0cm}{5mm}
\displaystyle v_h(0)=v_h(T)-\varphi_h,
\end{array} 
\right. \label{eqn_nonlocal}
\end{equation}
where $h=1/N$ and $v_h$ is a vector of approximations at spatial grid points 
\[ 
K_h=(hn_{1},hn_{2},\ldots,hn_{d}), ~~ 0\leq n_{j}\leq N, ~~ j=1,2,\ldots,d, 
\]
and $A_h$ is the second-order centered difference matrix operator\footnote{It is assumed that $A_h$ is symmetric and positive definite.} approximating the elliptic operator $-\nabla \cdot \left( a \nabla \right)$, incorporating also the homogeneous Dirichlet boundary conditions. By an additional time-stepping method, we can turn the continuous time solution $v_{h}(t)$ into a fully discrete one. Nevertheless, since problem \eqref{eqn_nonlocal} has a non-local condition in time, standard time stepping methods can not be exploited directly. Two natural ways of resolving the problem of non-locality in time is by 1) {\color{black} a direct discretization in time and space, which results in a linear system of extremely large size which is computationally expensive to solve, especially for two or three dimensional problems; \cite{Ashyralyev2012}, 2) developing a convergent shooting algorithm, see e.g., \cite{Geiger2015AdaptiveMS,CHEN2019226} for recent variants of it, to turn the non-local problem \eqref{eqn_nonlocal} into a purely initial-value problem, which is computationally and implementation-wise preferable in comparison to the former. In this artice, we develop low rank approximations based on the Arnoldi algorithm to efficiently transform the problem \eqref{eqn_nonlocal} into an initial-value problem, and accurately approximate the unknown pair $\{ p,u\}$,  in a much more efficient way than the fully discrete approach and the shooting method.}  

This article is structured as follows. In Section \ref{Sec_StandardMethods} we discuss two standard methods: a full discretization of the problem \eqref{eqn_nonlocal}, and a more efficient shooting algorithm, which is also proved to be convergent. Section \ref{Sec_Arnoldi_Approximation} contains the main contribution of this paper, where a low rank approximation algorithm based on Arnoldi is introduced, and the foundations of the ideas are solidified by a convergence analysis for the low rank approximations of the emerging operators in the algorithm. We conclude the paper by providing numerical experiments in Section \ref{Sec_Num_Results_TopSec}, to validate our theoretical findings.

\section{Preliminaries} \label{Sec_Preliminaries} 

To improve readability and to synchronize the mathematical notation throughout the paper, 1) we present the common notations that will be subsequently used in the sequel, 2) we present a few intermediate theoretical results in advance (to be used in the analysis later). 

\begin{itemize}
\item Let $B \in \mathbb{R}^{N^d\times N^{d}}$ be a matrix. We denote the set of eigenvalues of $B$ by $\lambda(B)$. The spectral radius of the matrix $B$ is denoted by $\rho_{B}$.
\item Throughout the paper, the letter $N$ represents the number of degrees of freedom in a spatial direction, while $M$ represents the degrees of freedom coming from a temporal discretization.
\item In order to avoid the confusion with the final time $T$ in equation \eqref{Eqn_u}, we denote the transpose of a matrix $B$ by $B^{*}$, even when the matrix $B$ is real-valued.
\item Superscripts are used to indicate the power of a matrix $B^{k}$, or the derivative of a function $f^{(k)}(t):=\frac{d^{k}f(t)}{dt^{k}}$, or as a means to represent the iteration number for a time-stepping method; see e.g. $v_h^{m}$ in equation \eqref{eqn_nonlocal_discrete}. The distinction between power of a quantity and the temporal iteration number is expected to be clear from the context. 

\item Properties of the matrix $A_h$:

\begin{itemize}
    \item[P1.] The matrix $A_h \in \mathbb{R}^{N^d \times N^d}$ is symmetric and positive definite. The eigenvalues $\lambda_j \in \lambda(A_h)$ are real and positive and ordered as
    \[ 
    0 < \lambda_1 \leq \lambda_2 \leq \ldots \leq \lambda_{N^d}, 
    \]
    and the corresponding eigenfunctions $\phi_j$ form an orthonormal basis for the finite dimensional space $l^2(K_h)$ equipped with the inner-product 
    \[ 
    \langle f_h,g_h \rangle = h^{d} \sum_{j=1}^{N^d} f_{h,j} g_{h,j}, 
    \]
    and the norm $\| f_h\|:= \sqrt{\langle f_h,f_h \rangle}$.
    \item[P2.] For sufficiently small $h$, we have $\lambda_1 \geq a_{\min} \frac{\pi^2}{d}$, up to an $O(h^2)$ discretization error.
    \item[P3.] The operator $A_h$ generates the contraction semigroup $e^{-TA_h}: l^2(K_h) \to l^2(K_h)$ which satisfies 
    \[  
    \left\|  e^{-TA_h} \right\| \leq e^{-\lambda_1 T}. 
    \]
\end{itemize}
{\bf Proof of P2. }  To see property $2$, let us consider the continuous operator $A = -\nabla \cdot \left(  a \nabla \right)$, which can be seen as the limit (as $h \to 0$) of the discrete operator $A_h$, and let $\{ \lambda_1,\phi_1 \}$ be the principal eigenvalue and eigen-function of $A$. Then it follows that  
\[ 
\int_{K} a(x) \nabla \phi_1(x) \cdot \nabla v(x) \; dx = \lambda_1 \int_{K} \phi_1(x) v(x) \; dx, \quad \forall v \in H_0^{1}(K). 
\]
Therefore, plugging $v = \phi_1$, we obtain
\[ 
\lambda_1 = \frac{\int_{K} a(x) \nabla \phi_1(x) \cdot \nabla \phi_1(x) dx}{\| \phi_1\|^{2}_{L^2(K)}} \geq a_{\min} \frac{\| \nabla \phi_1 \|^{2}_{L^2(K)}}{\| \phi_1\|^{2}_{L^2(K)}} \geq a_{\min} C_p(K)^{-2},
\]
where the last inequality follows from the Poincare's inequality, and $C_p(K)$ is the Poincare's constant for the domain $K$, which is bounded by $C_p(K) \leq \frac{diam(K)}{\pi} = \frac{\sqrt{d}}{\pi}$, see \cite{Payne1960}. Using this inequality, and the fact that $A_h \to A$, we obtain property $2$ for sufficiently small $h$. 

{\bf Proof of P3. }To see the last property, assume $g = \sum_{j=1}^{N^d} g_j \phi_j$, where $g_j =  \langle g,\phi_j \rangle$, and 
\begin{eqnarray*} 
\left\| e^{-TA_h} \right\|^2 & := & \sup_{\| g \| = 1} \left\| e^{-TA_h} g \right\|^2  
=  \sup_{\| g \| = 1} \|  \sum_{j=1}^{N^d} e^{-T \lambda_j } g_j \phi_j \|^2  \\ 
& = & \sup_{\| g \| = 1} h^{d} \sum_{j,\ell=1}^{N^{d}} e^{-T\lambda_j} e^{-T\lambda_{\ell}} g_j g_{\ell} \langle \phi_j,\phi_{\ell} \rangle. 
\end{eqnarray*}
Since the eigenfunctions are orthonormal, it follows that $\langle \phi_j,\phi_{\ell} \rangle = \delta_{i,j}$, where $\delta_{i,j}$ is the kronecker delta. Finally, using the discrete Parseval's equality, we obtain
\begin{eqnarray*}
\| e^{-TA_h} \|^2 & = & \sup_{\| g \| = 1}  h^{d} \sum_{j=1}^{N^{d}} e^{-2T\lambda_j}  g_j^2 \leq  \sup_{\| g \| = 1} e^{-2T\lambda_1} h^{d} \sum_{j=1}^{N^{d}} g_j^2 \\ 
& = & \sup_{\| g \| = 1}e^{-2T\lambda_1} \| g \|^2 = e^{-2T\lambda_1}.
\end{eqnarray*}
The result follows by taking the square roots of both sides.
\end{itemize}

\section{Standard numerical methods for solving \eqref{eqn_nonlocal}} \label{Sec_StandardMethods}

In this section, we present two strategies to tackle the problem of non-locality in time. First we present a direct approach, where the problem \eqref{eqn_nonlocal} is directly discretized in time, and we discuss potential numerical challenges with this approach. Second, we will present the shooting method, which is a more standard way of iteratively solving \eqref{eqn_nonlocal}. We conclude this section by a proof of convergence for the shooting method. 
 

\subsection{A Direct Approach} \label{SubSec_Direct_Approach}

Let $\tau=T/M$ and $t_{m}=m\tau$, $0\leq m \leq M$. A temporal discretization of the problem \eqref{eqn_nonlocal}, say using the  Crank-Nicholson scheme, yields
\begin{equation}
\left\{
\begin{array}{l}
\displaystyle \frac{v_h^{m+1}-v_h^{m}}{\tau}+\frac{A_{h}v_h^{m+1}+A_{h}v_h^{m}}{2}=f_{h}(t_{m+1/2}), ~ m=0,1,\ldots,M-1, \\ \rule{0cm}{6mm}
\displaystyle v_h^{0}=v_h^{M}-\varphi_h,
\end{array} 
\right. \label{eqn_nonlocal_discrete}
\end{equation}
where $v_{h}^{k}$ denotes the numerical approximations of $v_{h}(t_{k})$. Since $v_{h}^{0}$ is not available, \eqref{eqn_nonlocal_discrete} cannot be solved recursively. 
In fact, it is a system of $M(N-1)^{d}$ linear equations with $M(N-1)^{d}$ unknowns. Therefore, one needs to invert the square matrix of size $M(N-1)^{d}$ which requires a huge amount of computation and memory space; in particular, for multi-dimensional problems with sufficiently small step sizes $\tau$ and $h$. 

One can solve the nonlocal scheme \eqref{eqn_nonlocal_discrete} by using the locality of boundary conditions in one of the spatial variables. In fact, \eqref{eqn_nonlocal_discrete} can be written in the matrix form as follows:
\begin{equation} \label{eqn_matrix}
\left\{
\begin{array}{l}
\displaystyle CV_{n+1}+BV_{n}+CV_{n-1}=\phi_{n}, \quad 1\leq n\leq N-1,  \\ \rule{0cm}{5mm}
\displaystyle V_{0}=V_{N}=0,
\end{array}
\right. 
\end{equation}
where $B$ and $C$ are constant square matrices of size $M(N-1)^{d-1}$, $V_{n}$ and $\phi_{n}$ are column vectors with $M(N-1)^{d-1}$ entries. The solution of the matrix equation (\ref{eqn_matrix}) can be found by using the modified Gauss elimination method \cite{Samarskii_book} as follows:
\begin{equation*}
\left\{
\begin{array}{l}
\displaystyle V_{n}=\alpha_{n+1}V_{n+1}+\beta_{n+1}, \quad n=N-1,\ldots,2,1,  \\ \rule{0cm}{5mm}
\displaystyle V_{N-1}=0,
\end{array}
\right. 
\end{equation*}
where $\alpha_{n}$ are square matrices and $\beta_{n}$ are column vectors, calculated as
\begin{eqnarray*}
\alpha_{n+1} & = & -\left(B+C\alpha_{n}\right)^{-1}C, \\
\beta_{n+1} & = & \left(B+C\alpha_{n}\right)^{-1} \left(\phi_{n}-C\beta_{n}\right), ~~ n=1,2,\ldots,N-1. 
\end{eqnarray*}
Here, $\alpha_{1}$ is a zero matrix and $\beta_{1}$ is a zero vector.  

The elimination method described above requires the inversion of a square matrix of size $M(N-1)^{d-1}$ at every iteration, which significantly limits its application to one and two dimensional problems only. Due to this computational limitation, we argue that the use of fully discrete scheme should be avoided in practice, and that more efficient strategies for approximating \eqref{eqn_nonlocal} are needed. Finally, we note that other boundary conditions such as the Neumann conditions can also be treated with a slight adjustment but with no major conceptual change in the algorithm.


\subsection{Shooting Method} \label{SubSec_Shooting}

An alternative way to solve the non-local problem (\ref{eqn_nonlocal}) is by means of the shooting method, which is easy to implement and is favourable from a computational cost point of view. Let us denote by $v_h(t;\alpha)$ the solution of the initial value problem  
\begin{equation} \label{eqn_ivp}
\left\{
\begin{array}{l}
\displaystyle \frac{d}{dt}v_h(t;\alpha)+A_{h}v_h(t;\alpha)=f_h(t), \quad 0<t<T, \\ \rule{0cm}{5mm}
\displaystyle v_h(0;\alpha)=\alpha.
\end{array}
\right. 
\end{equation}
For $v_h(t;\alpha)$ to be a solution of (\ref{eqn_nonlocal}), the initial vector $\alpha$ must satisfy
\begin{equation} \label{eqn_FP_exact} 
\alpha=v_h(T;\alpha)-\varphi_h.
\end{equation}
The corresponding fixed-point iterations can be constructed as
\begin{equation} \label{eqn_FixedPoint_Iteration}
\alpha_{k+1}=v_h^{M}(\alpha_{k})-\varphi_h, \quad k=0,1,2,\ldots
\end{equation}
where $v_h^{M}(\alpha_{k})$ is the solution to the Crank-Nicholson approximation \eqref{eqn_nonlocal_discrete} of the problem \eqref{eqn_ivp} with initial data $\alpha_{k}$, which reads as 
\begin{eqnarray}
v_h^{m+1}(\alpha_{k}) & = & G(A_h;\tau) v_{h}^{m}(\alpha_k) + \tau \left(I + \frac{\tau A_h}{2}\right)^{-1}f_h(t_{m+1/2}), \label{eqn_Crank_Nicholson}  \\
v_h^{0}(\alpha_{k}) & = & \alpha_k, \quad m =0,1,2,\ldots,M-1, \quad \tau M = T, \nonumber
\end{eqnarray}
where $G(A_h;\tau):=\left(I + \frac{\tau A_h}{2}\right)^{-1} \left(I-\frac{\tau A_h}{2}\right)$.

\begin{rmk}
    {\color{black} Throughout the manuscript, we assume a Cranck-Nicholson discretization in time. The efficiency of the time stepping can be further improved by employing either explicit stabilized integrators with large stability regions, such as \cite{Vilmart_etal} or symplectic Runge Kutta methods such as \cite{LIU2021113133}. Nevertheless, note that employing different time-stepping methods will require completely different mathematical analysis, and the theories in this paper is valid only for the Cranck-Nicholson scheme. }
\end{rmk}

\begin{thm}
Let $\{\alpha_k\}$ be the sequence generated by fixed-point iterations \eqref{eqn_FixedPoint_Iteration}, and let $v_h(0)$ satisfy $v_h(0) = v_h(T;v_h(0)) - \varphi_h$, where $v_h(T;v_h(0))$ is the solution at time $T$ of problem \eqref{eqn_ivp} with initial data $v_h(0)$. Then it follows that 
\[
\|  \alpha_{k+1} - v_h(0) \| \leq \| G \|^{M} \| \alpha_{k} -  v_h(0) \| + C \tau^2 \left( 1+ \| G \|^{M} \right) \left(  1 - \| G \|^{M}\right)^{-1},
\] 
where $\| G \| = \left| \frac{1- \tau \lambda_1/2 }{1+ \tau \lambda_1/2 } \right| < 1$,  $\lambda_1$ is the principal eigenvalue of the matrix $A_h$, and $C$ is a constant independent of $k$ but may depend on $f_h$ and $v_h(0)$.  
\end{thm}

\begin{pf}
We re-write \eqref{eqn_Crank_Nicholson} in the form
\begin{eqnarray*}
v_h^{m+1}(\alpha_{k}) & = & G^{m+1} \alpha_k + \tau \sum_{j=0}^{m} G^{m-j} F_h(t_{j+1/2}), \\
F_h(t_{j+1/2}) & := & \left(I + \frac{\tau A_h}{2}\right)^{-1} f_h(t_{j+1/2}), \quad m =0,1,2,\ldots,M-1.
\end{eqnarray*}
Then the fixed-point iteration \eqref{eqn_FixedPoint_Iteration} can be written as 
\[ 
\alpha_{k+1}=G^{M} \alpha_k + \tau \sum_{j=0}^{M-1} G^{M-1-j} F_h(t_{j+1/2})-\varphi_h, \quad k=0,1,2,\ldots
\]
Moreover, let $\alpha_{\infty}$ be the limiting value of $\alpha_k$ satisfying\footnote{Indeed the fixed-point iteration is convergent since $\| G \| < 1$, and therefore the limiting value $\alpha_{\infty}$ exists.}  
\[ 
\alpha_{\infty}=G^{M} \alpha_{\infty} + \tau \sum_{j=0}^{M-1} G^{M-1-j} F_h(t_{j+1/2})-\varphi_h.
\]
Clearly 
\begin{equation} \label{Est_Intermediate_Thm1}
\| \alpha_{k+1} - \alpha_{\infty} \| \leq \| G\|^{M} \| \alpha_k  - \alpha_{\infty} \|.
\end{equation}
On the other hand,
\begin{eqnarray*}
\| \alpha_{\infty}  - v_h(0) \|_2  & = & \left\| v_h^{M}(\alpha_{\infty}) - \varphi_h - v_h(0) \right\| \\ 
& = & \left\| v_h^{M}(\alpha_{\infty}) - \varphi_h - \left( v_h(T;v_h(0))  - \varphi_h \right) \right\|_2  \\
& = & \left\| v_h^{M}(\alpha_{\infty}) -  v_h^{M}(v_h(0)) + v_h^{M}(v_h(0))  - v_h(T;v_h(0)) \right\|_2  \\
& \leq & \left\| v_h^{M}(\alpha_{\infty}) -  v_h^{M}(v_h(0)) \right\|_2  + \left\| v_h(T;v_h(0))  - v_h^{M}(v_h(0)) \right\|  \\
& \leq & \| G \|^{M} \| \alpha_{\infty} - v_h(0) \|  + C \tau^2.
\end{eqnarray*}
Therefore,
\[
\| \alpha_{\infty}  - v_h(0) \|  \leq C \tau^2 \left( 1- \| G \|^{M} \right)^{-1}. 
\]
Now using \eqref{Est_Intermediate_Thm1} together with the last estimate, we see that 
\begin{eqnarray*}
 \| \alpha_{k+1}  - v_h(0) \| & \leq & \| \alpha_{k+1} - \alpha_{\infty} \|  + \| \alpha_{\infty} - v_h(0) \| \\
 & \leq & \| G\|^{M} \| \alpha_k  - \alpha_{\infty} \| + \| \alpha_{\infty} - v_h(0) \| \\ 
 & \leq & \| G\|^{M} \| \alpha_k  - v_h(0)\| + \left(1 + \| G\|^{M}  \right) \| \alpha_{\infty} - v_h(0) \| \\
 & \leq & \| G\|^{M} \| \alpha_k  - v_h(0)\|  + C \tau^2 \left( 1- \| G \|^{M} \right)^{-1} \left(1 + \| G\|^{M}  \right).
\end{eqnarray*}
Together with the fact that $\| G\|  = \max_{j} \left| \frac{1- \tau \lambda_j/2 }{1+ \tau \lambda_j/2 } \right| = \left| \frac{1- \tau \lambda_1/2 }{1+ \tau \lambda_1/2 } \right| < 1 $, we conclude the proof.
$\square$
\end{pf}


\section{Arnoldi Approximation} \label{Sec_Arnoldi_Approximation}

Arnoldi approximation, \cite{Guntel2013,Higham_Book,Hochbruck1997}, relies on Krylov subspace methods and is typically used to approximate matrix functions $F(B)$ or  $F(B) b $, where $B \in \mathbb{C}^{N^d \times N^d}$ is a large matrix, and $b \in \mathbb{C}^{N^d}$ is a vector. The idea is to bypass the computation of large matrix functions and do the matrix computations on a lower dimensional subspace, which may result in tremendous computational gain. Consider, for example, the matrix exponential $F(B) = e^{-B}$, where $B \in \mathbb{R}^{N^d \times N^d}$ is a positive definite and symmetric matrix, which is also central to the goal of the present article. The starting point is a unitary transformation of the matrix $B$ in the form $H = Q^{*} B Q$, where  $Q \in \mathbb{R}^{N^d\times k}$ with $Q^{*}Q = I$, and $H \in  \mathbb{R}^{k\times k}$, where $k \ll N^d$. We can then approximate $F(B) b$ for $b\in \mathbb{R}^{N^d}$ as follows 
\[
F(B) b \approx  Q F(H) Q^{*} b.
\]
We summarize the properties of the matrix $H$, whenever the original matrix $B$ is symmetric and positive definite. We can make the following immediate observation: Assume that the matrix $Q$ has rank $k$. If the matrix $B$ is positive definite and symmetric, the matrix $H$ is also positive definite and symmetric. To prove the symmetric property, we write
\[
H^{*} = Q^{*} B^{*} Q =  Q^{*} B Q =H.
\]
To prove positiveness, assume $y\in \mathbb{R}^{k} \neq {\bf 0}$. Then 
\[
y^{*} H y =  y^{*} Q^{*} B Q y > 0.         
\]
The computational efficiency here originates from the fact that now the matrix exponential is computed for a matrix of much lower rank. The following theorem provides an error estimate for the difference $F(B) b - Q F(\mathcal{H}) Q^{*} b$. 

\begin{thm} \label{Thm_MatrixExponential} 
[Hochbruck,Lubich \cite{Hochbruck1997}] Let $B \in \mathbb{C}^{N^d\times N^d}$ be a Hermitian positive semi-definite matrix with eigenvalues in $[0,\rho]$. Moreover, let $H = Q^{\star} B Q$ be a unitary transformation of $B$ via an Arnoldi procedure with $H \in \mathbb{R}^{k \times k}$ and $Q \in \mathbb{R}^{N^d \times k}$. Then the following estimate holds
\begin{equation} \label{Arnoldi_Exponential_Estimate}
\|  e^{-B} b  - Q e^{-H} Q^{\star} b \|_2 \leq \| b \|_2 
\left\{
\begin{array}{ll}
\displaystyle 10 e^{-\frac{4 k^2}{5 \rho}}, & \sqrt{\rho} \leq k \leq \rho/2, \\ \rule{0cm}{7mm}  
\displaystyle \frac{40}{\rho} e^{-\frac{\rho}{4}} \left( \frac{e \rho}{4 k} \right)^{k}, & k \geq \rho/2.
\end{array}
\right.
\end{equation}
\end{thm}

Theorem \ref{Thm_MatrixExponential} results in a particular computational advantage when applied to the second order difference operator $A_h$. The spectral radius $\rho_{A_h}$ of the matrix $A_h$ scales as $h^{-2}$. This will then imply that in order to obtain an exponential accuracy of the form $e^{-C}$ (for some positive constant $C$), it suffices to choose $k = \sqrt{5C/4} h^{-1}$. In other words, regardless of the dimension of the problem, we can compute $e^{-A_h} b$, where $b\in \mathbb{R}^{N^d}$ at the cost of computing the matrix exponential of a one dimensional problem. 

In what follows, we will present two different ways of exploiting the Arnoldi method to gain computational efficiency in comparison to the direct approach from section \ref{SubSec_Direct_Approach}, and the shooting method from section \ref{SubSec_Shooting}. Both approaches rely on relating the solution $v(T)$ of the parabolic PDE \eqref{Eqn_v} at time $T$ to the initial value $v(0)$ and then approximating the emerging operators using a low rank Arnoldi approximation. This will then result in significantly lower computational cost while retaining the same accuracy of the direct discretization and the shooting method. For both approaches, we also include separate convergence analysis, which shows that the emerging matrix functions can be approximated accurately at the cost of matrix function computations of a one-dimensional problem.     

\begin{rmk}
    Note that matrix exponentials of the form $e^{-B}$ can also be computed by a standard eigenfunction expansion. This is, however, computationally very expensive and must be avoided in computations, see \cite{Abdulle_Arjmand_Paganoni_22} for a full error analysis. 
\end{rmk}


\subsection{A hybrid Shooting-Arnoldi approximation} \label{Sec_Approach1} 

{\bf Step 1.} The starting point is to write the solution $v_h(T)$ of \eqref{eqn_nonlocal} in terms of the initial value and the right hand side


\begin{equation} \label{Eqn_Duhammel}
v_h(T) = e^{-TA_h} v_h(0) + \int_{0}^{T} e^{-(T-s) A_h } f_h(s) ds.
\end{equation}

{\bf Step 2.} Approximate the operators in \eqref{Eqn_Duhammel} using the Arnoldi algorithm, and use the fixed point iteration to find $v_h(0)$ such that $v_h(0)  = v_h(T) - \varphi_h$ is met, up to a desired tolerance. In other words, solve for 

\begin{equation} \label{Eqn_FixedPoint}
\alpha_{n} =  Q e^{-TH} Q^{*} \alpha_{n-1}  + \int_{0}^{T} Q e^{-(T-s) H} Q^{*} f_{h}(s) ds  - \varphi_h,     
\end{equation}
for sufficiently large $n$ such that $\alpha_n \approx v_h(0)$.

{\bf Step 3.} Solve \eqref{eqn_nonlocal} as an initial value problem, where the initial data comes from Step 2.

\subsubsection{Analysis}

Here, we aim at establishing an error bound for the difference between $\alpha_n$, defined by \eqref{Eqn_FixedPoint}, and the true value of the initial data. The main result is stated in Lemma \ref{Lem_Approcah_One}.

\begin{lem}\label{Lem_Fixed_Point} 
Suppose $A_h \in \mathbb{R}^{N^d \times N^d}$ is the second order difference matrix in \eqref{eqn_nonlocal}. Let $v_h(0)$ be the vector satisfying the relations \eqref{Eqn_Duhammel} and $v_h(0) = v_h(T) - \varphi_h$ exactly, and $\{\tilde{\alpha}_n\}$ be the sequence generated by the fixed point iterations
\begin{equation} \label{Eqn_FixedPointwithAh}
 \tilde{\alpha}_{n} = e^{-TA_h} \tilde{\alpha}_{n-1}  + \int_{0}^{T} e^{-(T-s)A_h} f_h(s) ds  - \varphi_h, \quad n=1,2,\ldots.
\end{equation}
Then it follows that 
\[
 \| \tilde{\alpha}_{n} - v_h(0) \| \leq e^{-Tn \lambda_1} \| \tilde{\alpha}_{0} - v_h(0) \|, \quad n=1,2,\ldots,  
\]
where $\lambda_1 > a_{\min} \frac{\pi^2}{d}$ is the smallest eigenvalue of the operator $A_h$. 
\end{lem}

\begin{pf} Using the relation $v_h(0) = v_h(T) - \varphi_h$ together with \eqref{Eqn_Duhammel} we see that
\[
v_h(0) = e^{-TA_h} v_h(0)  + \int_{0}^{T} e^{-(T-s)A_h} f_{h}(s) ds  - \varphi_h.
\]
Using the last equation and \eqref{Eqn_FixedPointwithAh}, we immediately see that 
\begin{eqnarray*}
\|  \tilde{\alpha}_n - v_h(0) \| & = & \| e^{-TA_h} \left( \tilde{\alpha}_{n-1} -  v_h(0) \right) \| \\ 
& \leq & \| e^{-TA_h} \| \| \tilde{\alpha}_{n-1} -  v_h(0)  \| \leq  e^{-T\lambda_1}  \| \tilde{\alpha}_{n-1} -  v_h(0)  \|, 
\end{eqnarray*}
where property $3$ in Section \ref{Sec_Preliminaries} was used for establishing the last inequality. $\square$
\end{pf}

\begin{lem} \label{Lem_Approcah_One} 
Let $A_h \in \mathbb{R}^{N^d \times N^d}$ be the same matrix as in Lemma \ref{Lem_Fixed_Point} and $0<\varepsilon<T$. Suppose $Q \in \mathbb{R}^{N^d \times k}$ and $H \in \mathbb{R}^{k \times k}$ with 
\begin{equation} \label{AssumptionOnk}
\sqrt{\frac{5T\rho_{A_h}}{4} \ln{\frac{10}{1-e^{-T\lambda_1}}}}<k<\frac{\varepsilon\rho_{A_h}}{2}
\end{equation}
are the corresponding matrices coming from an Arnoldi procedure applied to the matrix $A_h$. Moreover, let $v_h(0)$ be the vector satisfying the relations \eqref{Eqn_Duhammel} and $v_h(0) = v_h(T) - \varphi_h$ exactly, and $\alpha_n$ be the solution of the fixed point iteration \eqref{Eqn_FixedPoint}, with $\alpha_0 = \tilde{\alpha}_0$, where $\tilde{\alpha}_0$ is the initial guess for the iteration \eqref{Eqn_FixedPointwithAh}. Then it follows that 
\begin{eqnarray}\label{Main_Estimate_Full_FixedPoint}
\| \alpha_{n} - v_h(0) \| & \leq & C^n \| \alpha_{0} - v_h(0) \|  \nonumber \\ 
& + & 10 e^{-\frac{4 k^2}{5T\rho_{A_h}}} \left(\|v_h(0)\| + T\sup_{0\leq s \leq T} \|f_h(s)\|  \right) \frac{1-C^{n}}{1-C} \nonumber \\
& + &  \varepsilon \left(1+ \| Q\| \| Q^{*}\| \right) \sup_{0\leq s \leq T} \|f_h(s)\|  \frac{1-C^{n}}{1-C},
\end{eqnarray}
where $C = e^{-T \lambda_1} + 10 e^{-\frac{4 k^2}{5T\rho_{A_h}}}<1$. 
\end{lem}

\begin{pf} Let $\tilde{\alpha}_n$ and $v_h(0)$ be defined in the same way as in the proof of Lemma \ref{Lem_Fixed_Point}. We start by splitting the error $\alpha_{n} - v_h(0)$ into two parts as follows: 
\begin{equation} \label{TE_Estimate_Full_FixedPoint}
\| \alpha_{n} - v_h(0) \| \leq \| \alpha_{n} - \tilde{\alpha}_n \| + \| \tilde{\alpha}_n - v_h(0) \|.
\end{equation}
An upper bound for the second term in the right hand side was already established in Lemma \ref{Lem_Fixed_Point}. To bound the first term in the right hand side, we write
\begin{eqnarray*}
\tilde{\alpha}_{n} - \alpha_n & = & e^{-TA_h} \tilde{\alpha}_{n-1}  + \int_{0}^{T} e^{-(T-s)A_h} f_h(s) ds  - \varphi_h \\ 
& - & \left( Q e^{-TH} Q^{*} \alpha_{n-1} + \int_{0}^{T} Q e^{-(T-s) H} Q^{*} f_h(s) - \varphi_h \right) \\ 
& = & e^{-TA_h} \left( \tilde{\alpha}_{n-1} - \alpha_{n-1}\right) - \left( e^{-TA_h} - Q e^{-TH} Q^* \right) \left( \tilde{\alpha}_{n-1} - \alpha_{n-1} \right) \\ 
& + & \left( e^{-TA_h} - Q e^{-TH} Q^* \right) \left( \tilde{\alpha}_{n-1}  - v_h(0)\right) \\ 
& + & \left( e^{-TA_h}  - Q e^{-TH}Q^*\right) v_h(0) \\ 
& + & \int_{0}^{T} \left( e^{-(T-s) A_h} - Q e^{-(T-s)H} Q^{*} \right) f_h(s) ds.
\end{eqnarray*}
Then,
\begin{eqnarray*}
\| \tilde{\alpha}_{n} - \alpha_n \| & \leq & \left( \| e^{-TA_h} \| + \|e^{-TA_h} - Q e^{-TH} Q^* \| \right) \| \tilde{\alpha}_{n-1} - \alpha_{n-1} \| \\
& + & \|e^{-TA_h} - Q e^{-TH} Q^* \| \| \tilde{\alpha}_{n-1} - v_h(0) \|  \\
& + & \|e^{-TA_h} - Q e^{-TH} Q^* \| \| v_h(0)\| \\ 
& + & \int_0^{T-\varepsilon} \|e^{-(T-s) A_h} - Q e^{-(T-s)H} Q^{*} \| \| f_h(s) \| ds \\ 
& + & \int_{T-\varepsilon}^{T} \|e^{-(T-s) A_h} - Q e^{-(T-s)H} Q^{*} \| \| f_h(s) \| ds.
\end{eqnarray*}
Now, using the fact that $\| e^{-TA_h} \|  \leq e^{-T \lambda_1}$ (see property $3$ in Section \ref{Sec_Preliminaries}), as well as $\| e^{-(T-s)A_h} - Q e^{-(T-s)H} Q^* \| \leq 10 e^{-\frac{4 k^2}{5(T-s)\rho_{A_h}}}, 0 \leq s \leq T-\varepsilon$ (see Theorem \ref{Thm_MatrixExponential}), we obtain
\begin{eqnarray*}
\| \tilde{\alpha}_{n} - \alpha_n \| & \leq & \left(  e^{-T \lambda_1}  + 10 e^{-\frac{4 k^2}{5T\rho_{A_h}}} \right) \| \tilde{\alpha}_{n-1} - \alpha_{n-1} \| \\
& + & 10 e^{-\frac{4 k^2}{5T\rho_{A_h}}} \| \tilde{\alpha}_{n-1} - v_h(0) \| + 10 e^{-\frac{4 k^2}{5T\rho_{A_h}}} \| v_h(0)\|  \\ 
& + &   10 e^{-\frac{4 k^2}{5T\rho_{A_h}}} (T-\varepsilon)\sup_{0\leq s \leq T-\varepsilon} \| f_h(s) \| \\
& + & \varepsilon \left(1+ \| Q\| \| Q^{*}\| \right) \sup_{T-\varepsilon\leq s \leq T} \| f_h(s) \|.
\end{eqnarray*}
We note that the assumption \eqref{AssumptionOnk} on $k$ implies that the amplification factor $C = e^{-T \lambda_1} + 10 e^{-\frac{4 k^2}{5T\rho_{A_h}}}<1$. 
Moreover, using Lemma $1$ for the term $\| \tilde{\alpha}_{n-1} - v_h(0) \|$, we get
\begin{eqnarray*}
\| \tilde{\alpha}_{n} - \alpha_n \| & \leq & C \| \tilde{\alpha}_{n-1} - \alpha_{n-1} \|  + 10 e^{-\frac{4 k^2}{5T\rho_{A_h}}} e^{-T(n-1) \lambda_1}\| \tilde{\alpha}_{0} - v_h(0) \|  \\ 
& + &  10 e^{-\frac{4 k^2}{5T\rho_{A_h}}} \left(\| v_h(0)\| + T\sup_{0\leq s \leq T} \| f_h(s) \| \right) \\
& + & \varepsilon \left(1+ \| Q\| \| Q^{*}\| \right) \sup_{0\leq s \leq T} \| f_h(s) \|.
\end{eqnarray*}
Then, by a trivial induction we have 
\begin{eqnarray*}
\| \tilde{\alpha}_{n} - \alpha_n \| & \leq & C^{n} \| \tilde{\alpha}_{0} - \alpha_{0} \|  + 10 e^{-\frac{4 k^2}{5T\rho_{A_h}}} \| \tilde{\alpha}_{0} - v_h(0) \| \sum\limits_{\ell=0}^{n-1} C^{\ell} e^{-T(n-1-\ell) \lambda_1}\\ 
& + &  10 e^{-\frac{4 k^2}{5T\rho_{A_h}}} \left(\| v_h(0)\| + T\sup_{0\leq s \leq T} \| f_h(s) \| \right) \sum\limits_{\ell=0}^{n-1} C^{\ell} \\
& + & \varepsilon \left(1+ \| Q\| \| Q^{*}\| \right) \sup_{0\leq s \leq T} \| f_h(s) \| \sum\limits_{\ell=0}^{n-1} C^{\ell}.
\end{eqnarray*}
The final estimate follows then from \eqref{TE_Estimate_Full_FixedPoint} under the assumption $\tilde{\alpha}_0 = \alpha_0$. 
\end{pf}

\begin{rmk}
    Note that the estimate \eqref{Main_Estimate_Full_FixedPoint} in Lemma \ref{Lem_Approcah_One} is sub-optimal due to the presence of an additional $O(\varepsilon)$ term in the upper bound. Indeed, in numerical simulations we only observe the exponentially decaying part of the error, see the numerical results in subsection \ref{Sec_MatrixFunction_Approx}.
\end{rmk}

\subsection{An approach based on a pure Arnoldi approximation} \label{Sec_Approach2} 

{\bf Step 1.} Unlike the previous approach, instead of using iterations to solve for $v_h(0)$, one can also formulate a direct equation by exploiting \eqref{Eqn_Duhammel}, together with $v_h(0)  = v_h(T) - \varphi_h$ to represent the initial data $v_h(0)$ as follows

\begin{equation} \label{Eqn_v0_intermsofvT}
v_h(0) = (I - e^{-TA_h})^{-1} \left(  -\varphi_h  +  \int_{0}^{T} e^{-(T-s) A_h } f_h(s) ds \right).
\end{equation}

{\bf Step 2.} Approximate the matrix operators $(I - e^{-TA_h})^{-1}$, and $e^{-A_h (T-s)}$ using the Arnoldi algorithm, and solve for $v_h(0)$.

{\bf Step 3.} Solve \eqref{eqn_nonlocal} as an initial value problem, where the initial data comes from Step 2. 

\subsubsection{Analysis}

In this approach, in addition to approximations of exponential functions, we also need to approximate the operator $(I - e^{-TA_h})^{-1}$. In what follows, we aim at establishing error bounds for an Arnoldi approximation for this operator. The main result is Theorem \ref{Thm_Estimate_1}. However, we will need the following lemma prior to proving Theorem \ref{Thm_Estimate_1}.




\begin{lem} [See Chapter 11, \cite{GolubVanLoan13}] \label{Lemma_Matrix_Power_Series}
If a function $f$ has a power series representation
\[
f(z)  =  \sum_{k=0}^{\infty} c_k z^k
\]
on an open disk containing the eigenvalues $\lambda(B)$ of a matrix $B \in \mathbb{C}^{N^d\times N^d}$, then 
\[
f(B)  =  \sum_{k=0}^{\infty} c_k B^k,
\]
and, 
\begin{equation} \label{Estimate_MatrixFunction_Trunc}
\|  f(B)  - \sum_{j=0}^{k-1} c_j B^{j} \| \leq \frac{N^d}{k!} \max_{0\leq s \leq 1} \| B^k f^{(k)}(Bs) \|.
\end{equation}
\end{lem}
Now suppose $b(z)  = (1-  e^{-z})^{-1}$. We are interested in deriving an error estimate for 
\[
\|  b(T A_h) g  - Q b(T H) Q^{*} g \|.
\]

\begin{thm} \label{Thm_Estimate_1} 
Let $A_h \in \mathbb{R}^{N^d \times N^d}$ be as in Lemma \ref{Lem_Approcah_One}, and $Q \in \mathbb{R}^{N^d \times k}$ and $H\in \mathbb{R}^{k\times k}$ with $\mu_1  = \min \lambda(H)$, be the matrices coming from the Arnoldi procedure. Then the following estimate holds:  
\begin{eqnarray*}
\|  b(T A_h) g  & - & Q b(T H) Q^{*} g \|  \leq  \left(N^{d} e^{-T \lambda_1 {m}} \frac{1}{(1 - e^{- T\lambda_1})^{m+1}} \right. \\ 
 & + & \left. C k e^{-T \mu_1 {m}} \frac{1}{(1 - e^{- T\mu_1})^{m+1}} 
 +  10 (m-1)  e^{-\frac{4 k^{2-\beta}}{5 \rho}} \right) \| g \|,
\end{eqnarray*}
where $m = \lfloor k^{\beta} \rfloor$, $0<\beta<1$, and the constant $C = \| Q\| \| Q^{*}\|$.
\end{thm} 

\begin{pf}
Let $f(z) = (1-z)^{-1}$. Then for $|z| < 1$, we can write 
\[
f(z)  = \sum_{j=0}^{\infty} z^j  = f_k(z)  + \sum_{j=k}^{\infty} z^j, \quad |z|<1,
\]
where $f_k(z)  = \sum\limits_{j=0}^{k-1} z^j$. Moreover, let $B = e^{-TA_h}$, then clearly $\lambda(B) \in [e^{-T\rho_{A_h}}, e^{-T\lambda_1}]$, and by Lemma \ref{Lemma_Matrix_Power_Series} it follows that 
\[
b(TA_h)  = f(B) = f_k(B)  + \sum\limits_{j=k}^{\infty} B^j,
\]
where $f_{k}(B)$ is given by 
\[
f_{k}(B)  = \sum\limits_{j=0}^{k-1} B^{j}.
\]
Note that the spectral radius of $B$ is strictly less than $1$. Now, assume $0< \beta < 1$ and $m=\lfloor k^{\beta} \rfloor$, and consider the  following decomposition of the overall error:
\[
\|  f(B) g  - Q f(e^{-TH}) Q^{*} g \| \leq \| f(B) g - f_{m}(B) g \|  + \| f_{m}(B) g  - Q f(e^{-TH}) Q^{*} g  \|.
\]
To bound the first term, we use inequality \eqref{Estimate_MatrixFunction_Trunc} in Lemma \ref{Lemma_Matrix_Power_Series} 
\begin{eqnarray*}
\| f(B) g - f_{m}(B) g \| & \leq & \| f(B) - f_{m}(B) \| \| g \| \\
& \leq & \frac{N^d}{m!} \max_{0\leq s \leq 1} \| B^{m} f^{(m)}(Bs) \| \| g \| \\ 
& \leq & \frac{N^d}{m!} \| B^{m}  \| \max_{0\leq s \leq 1} \| f^{(m)}(Bs) \| \| g \| \\ 
& \leq & \frac{N^d}{m!}e^{-T \lambda_1 {m}} \max_{0\leq s \leq 1} \| \frac{m!}{(I - B s)^{m+1}} \| \| g \| \\ 
& \leq & N^d e^{-T \lambda_1 {m}} \frac{1}{(1 - e^{- T\lambda_1})^{m+1}} \| g \|.
\end{eqnarray*}
For the second term, we consider the decomposition
\begin{eqnarray*}
\left\| f_{m}(B) g  - Q f\left(e^{-TH}\right) Q^{*} g  \right\|  & \leq & \left\| f_{m}(B) g  - Q f_{m}\left(e^{-TH}\right) Q^{*} g  \right\| \\
& + & \left\| Q \left( f\left(e^{-TH}\right)-f_{m}\left(e^{-TH}\right)\right) Q^{*} g  \right\|.
\end{eqnarray*}
The first term in the right hand side can be bounded as follows
\begin{eqnarray*}
\left\| f_{m}(B) g  - Q f_{m}\left(e^{-TH}\right) Q^{*} g  \right\| & \leq &
\sum_{\ell=1}^{m-1} 10 \| g \| e^{-\frac{4 k^2}{5 \ell \rho}} \\
& \leq & 10 (m-1) \| g \| e^{-\frac{4 k^2}{5 (m-1) \rho}} \\ 
& \leq & 10 (m-1) \| g \| e^{-\frac{4 k^{2-\beta}}{5 \rho}}.
\end{eqnarray*}
The second term, on the other hand, is bounded in a similar way as $\| f(B) g - f_{m}(B) g \|$ with an additional constant accounting for the $l^2$ norms of $Q$ and $Q^*$, i.e., 
\[
\left\| Q \left( f\left(e^{-TH}\right)-f_{m}\left(e^{-TH}\right)\right) Q^{*} g  \right\| \leq C k e^{-T \mu_1 m} \frac{1}{\left(1 - e^{- T\mu_1}\right)^{m}} \|g\|. 
\]

\end{pf}

\begin{rmk}
Choosing $0< \beta <1$ sufficiently small and $k = O(N)$, Theorem \ref{Thm_Estimate_1} results in exponentially decaying error bounds in terms of $N$, while the cost would be comparable to computing the matrix functions of sizes $k\times k$.  
\end{rmk}

\subsection{Numerical integration of matrix exponentials}

In the analysis provided in Sections \ref{Sec_Approach1} and \ref{Sec_Approach2}, it was assumed that the integral $\int_{0}^{T} e^{-(T-t)A_h} f_h(t) \; dt$ and it's approximation $\int_{0}^{T} Q e^{-(T-t)H} Q^{*} f_h(t) \; ds$ are given exactly. In practice these integrals need to be approximated by a quadrature rule. Here, we draw the attention of the reader to a subtlety in relation with the quadrature rule. For simplicity, we consider only the integral $\int_{0}^{T} e^{-(T-s)A_h} f_h(s) \; ds$.  A quadrature rule based on the standard midpoint rule yields
 \[
\int_0^T e^{-(T-t) A_h } f_h(t) dt = \tau \sum_{k=1}^{M} e^{-(T-t_{k-1/2}) A_h } f_h(t_{k-1/2}) + O(\tau^2).
 \]
 The error bound for such an approximation includes the second derivative of the integrand $e^{-(T-t) A_h } f_h(t)$, which would scale as $\| A^2\| \approx (1/h)^4 \approx N^{4}$. Therefore, a direct application of the midpoint rule will be very costly since very small time steps will be needed to achieve reasonable error tolerances. The problematic part with the error bound for the standard midpoint rule is that the exponential term in the integrand is also approximated at the midpoint. We will now, instead, derive a variant of the midpoint rule which bypasses the approximation of the exponential part of the integrand. This would then result in a much smaller pre-factor in the error bound as we will see in a nutshell. To achieve this, we start by
\begin{eqnarray*}
\int_0^T e^{-(T-t)A_h } f_h(t)  dt & = & \sum_{k=1}^{M} \int_{t_{k-1}}^{t_k} e^{-(T-t)A_h} f_h(t) dt \\
& \approx & \sum_{k=1}^{M} \int_{t_{k-1}}^{t_k} e^{-(T-t) A_h}\; dt  f_h(t_{k-1/2})  
\end{eqnarray*}
Moreover, 
\[
\int_{t_{k-1}}^{t_k} e^{-(T-t)A_h} dt = e^{-TA_h} \int_{t_{k-1}}^{t_k} e^{t A_h} dt = e^{-TA_h} \left(  e^{t_{k} A_h} -  e^{t_{k-1} A_h}  \right) A_h^{-1}.   
\]
Therefore,
\[
\int_0^T e^{-(T-t)A_h } f_h(t) dt = \sum_{k=1}^{M} \left( e^{-(T- t_k)A_h } - e^{-(T-t_{k-1}) A_h } \right) A_h^{-1} f_h(t_{k-1/2}) + O(\tau). 
\]
The procedure above results in a first order method but the pre-factor in the upper bound is indepedent of $A$, and includes only the first derivative $\partial_t f_h(t)$. In our simulations below, we also use a Richardson extrapolation to achieve a globally second-order method in time, which proves to be much more efficient than a standard midpoint rule.

\section{Numerical Results} \label{Sec_Num_Results_TopSec}

In this section, we provide numerical examples to validate the theoretical findings. In subsection \ref{Sec_MatrixFunction_Approx}, we provide a numerical example to demonstrate the exponential convergence for the approximation of the matrix exponential $e^{-TA_h}$ as well as $(I - e^{-TA_h})^{-1}$. In subsections \ref{Sec_Numerical_Results_up} we provide numerical examples in one, two, and three dimensions to demonstrate the advantage of using the full Arnoldi approximation or the hybrid approach in comparison to the standard shooting method.

\subsection{Matrix function approximations-Numerical results} \label{Sec_MatrixFunction_Approx}

\begin{figure}[h]
\centering
  \includegraphics[width=10cm]{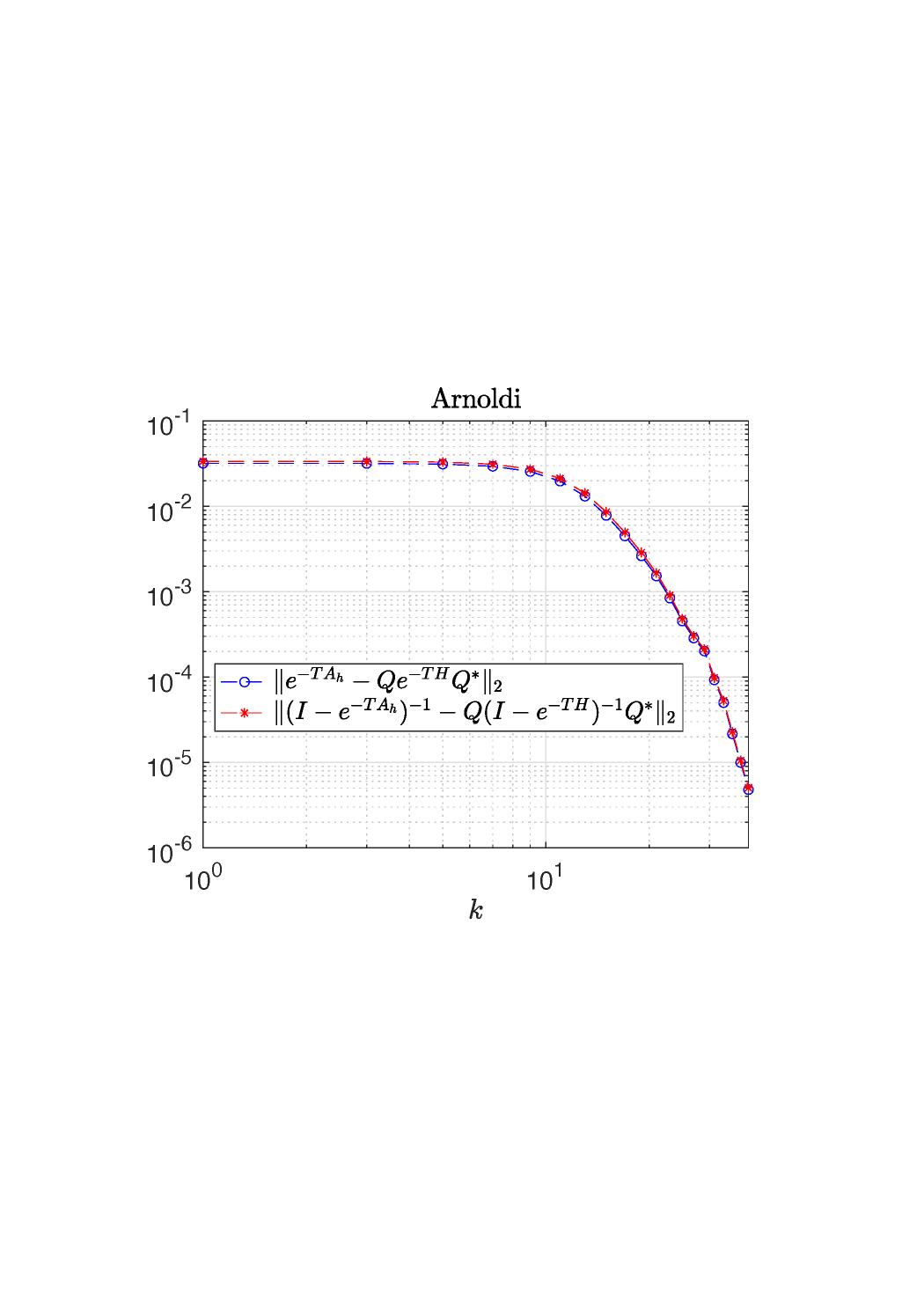}
\vspace{-4cm}
  \caption{Approximation error for the matrix functions computations via Arnoldi.}
  \label{fig:Arnoldi_Conv}
\end{figure}

In this section, we provide numerical evidence for exponentially decaying errors originating from the Arnoldi approximation of the operators $e^{-T A_h}$, and $(I - e^{-TA_h})^{-1}$. The size of the matrix $A_h \in  \mathbb{R}^{N^2 \times N^2}$ is $1600 \times 1600$. Figure \ref{fig:Arnoldi_Conv} demonstrates the exponential decay of the error corroborating the results of Theorems \ref{Thm_MatrixExponential} and \ref{Thm_Estimate_1}. In particular, we observe accuracies down to $10^{-5}$ tolerances for reasonably low rank approximation; i.e., $k\approx 40$. Note that, in this simulation the value of the parameter $T$ is chosen to be $T = 0.1$. 

\subsection{Approximation of $u$ and $p$} \label{Sec_Numerical_Results_up}

\begin{figure}
\centering
\begin{subfigure}{.5\textwidth}
  \centering
  \includegraphics[width=.8\linewidth]{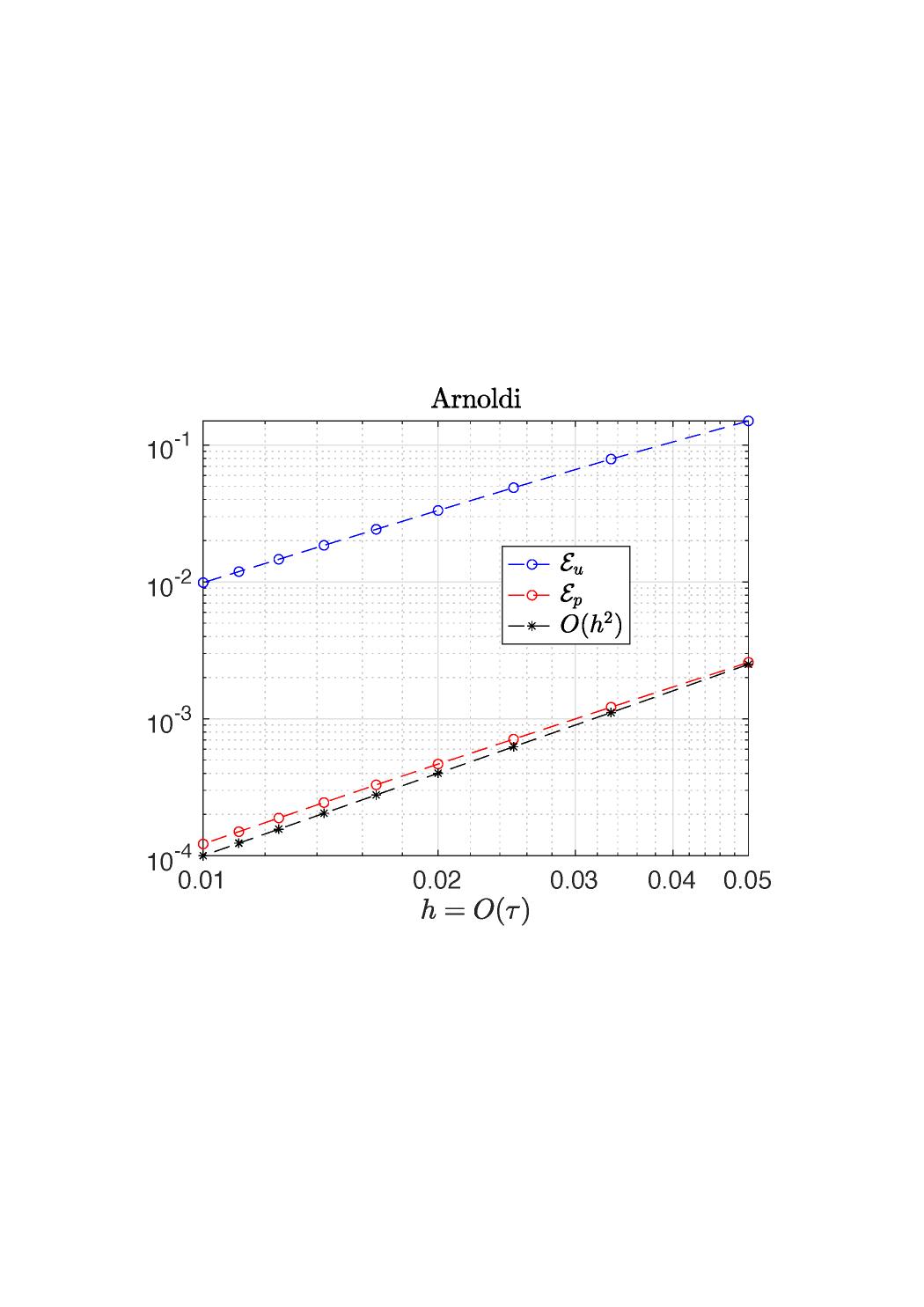}
  \caption{Convergence of $u$ and $p$ in $d=1$}
  \label{fig:OneD_Ex_1}
\end{subfigure}%
\begin{subfigure}{.5\textwidth}
  \centering
  \includegraphics[width=.8\linewidth]{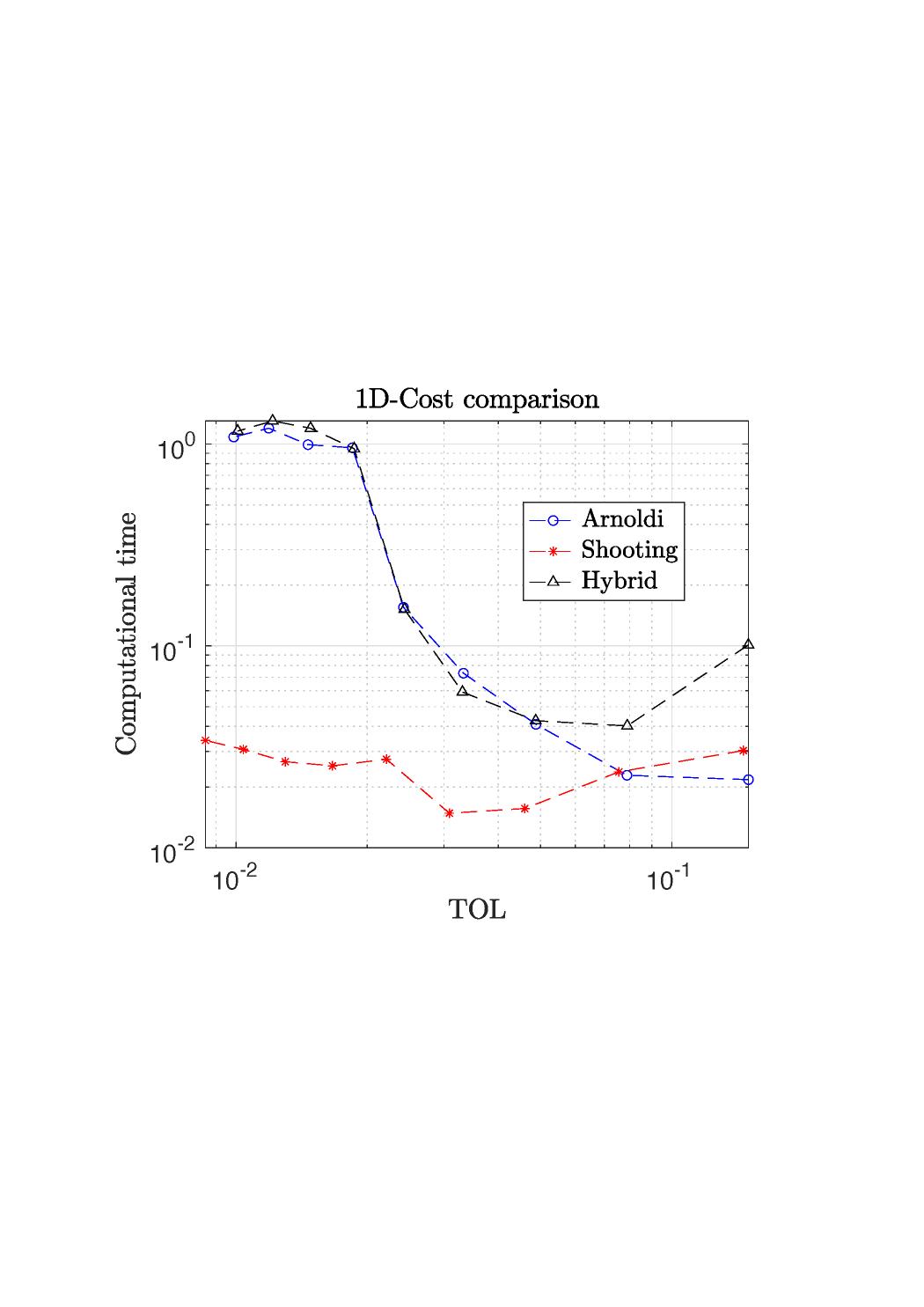}
  \caption{Cost comparison in $d=1$}
  \label{fig:OneD_Cost}
\end{subfigure}
\caption{Numerical results in one dimension.}
\label{fig:OneDResults}
\end{figure}

\begin{figure}
\centering
\begin{subfigure}{.5\textwidth}
  \centering
  \includegraphics[width=.8\linewidth]{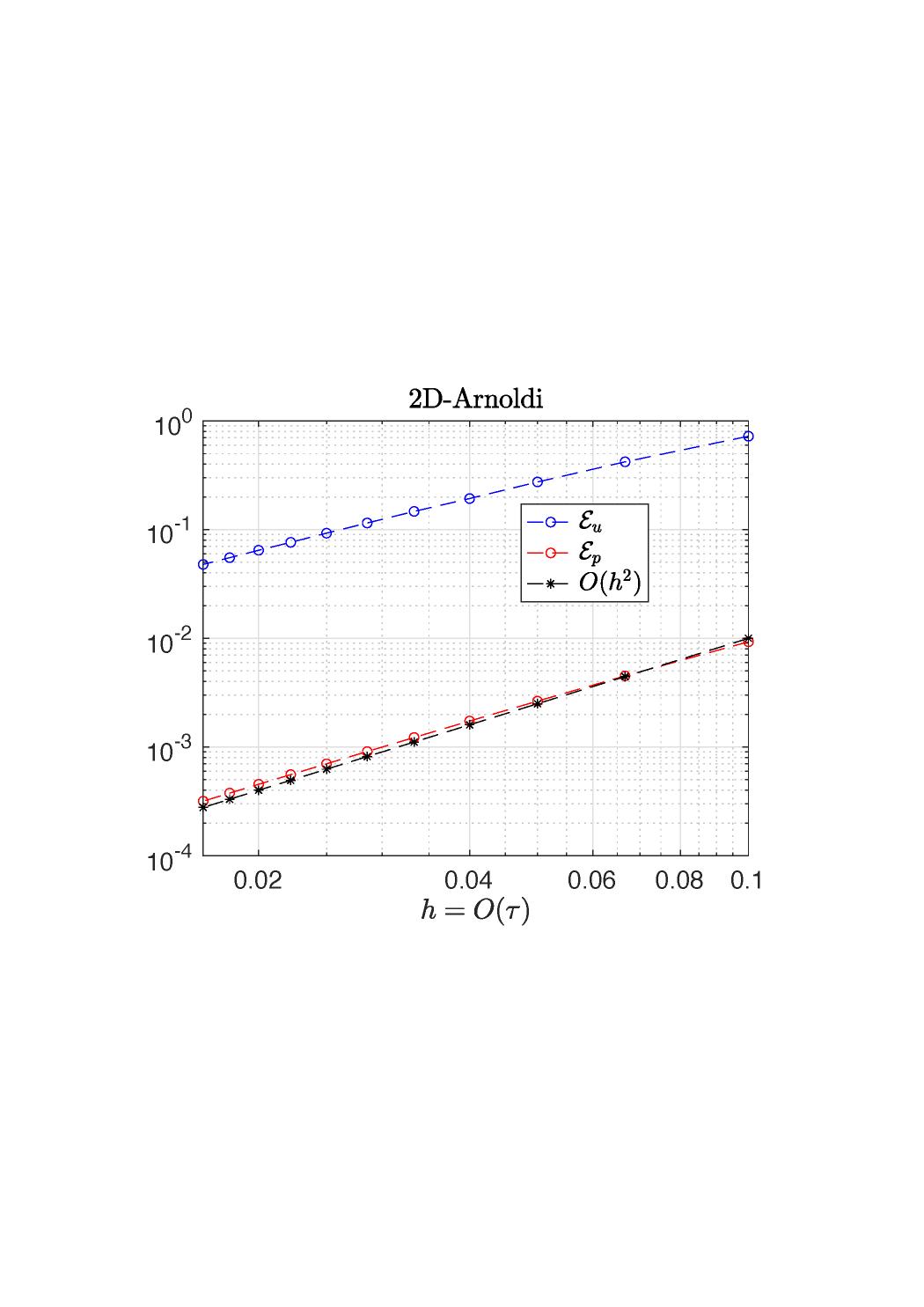}
  \caption{Convergence of $u$ and $p$ in $d=2$}
  \label{fig:TwoD_Ex_1}
\end{subfigure}%
\begin{subfigure}{.5\textwidth}
  \centering
  \includegraphics[width=.8\linewidth]{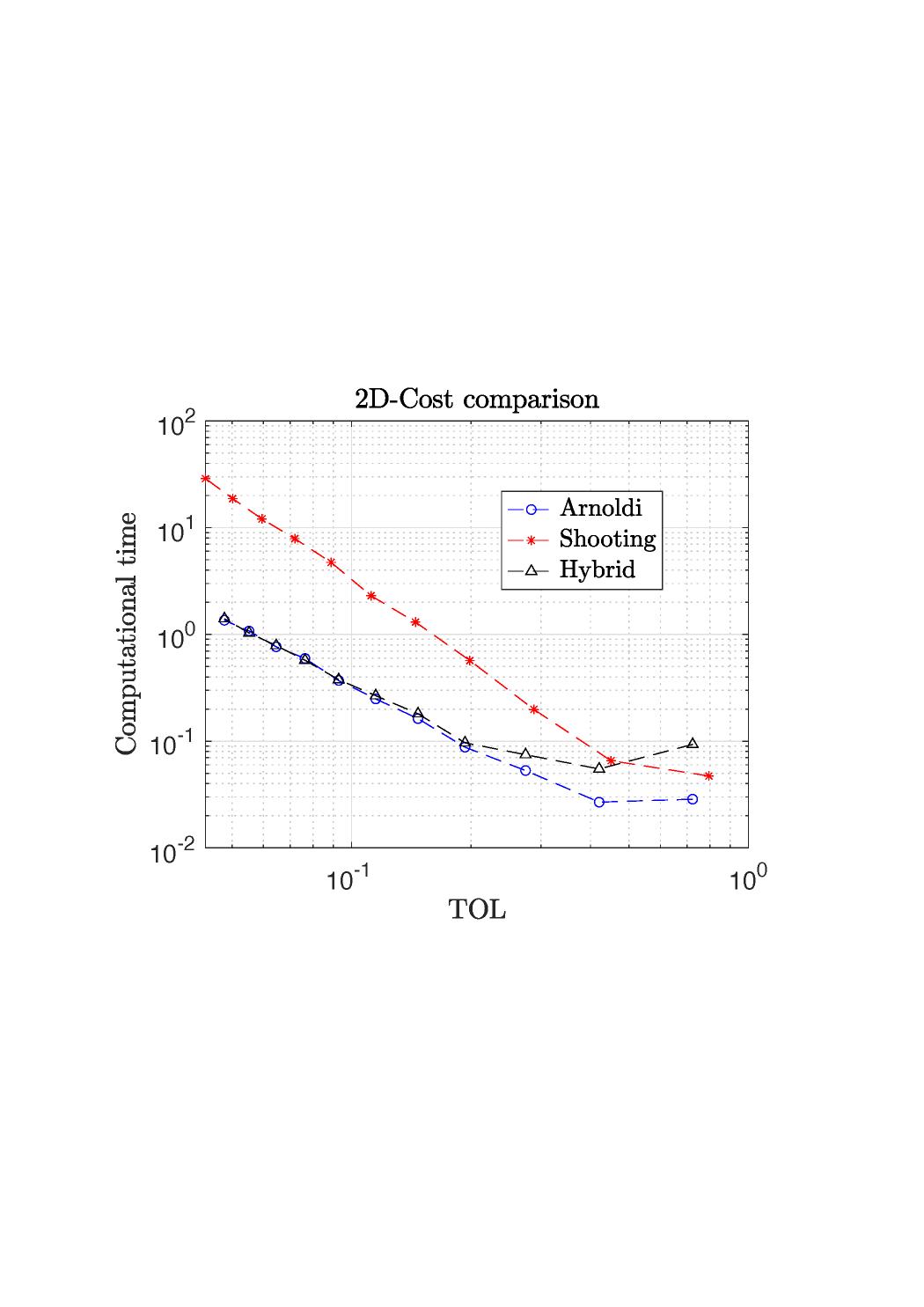}
  \caption{Cost comparison in $d=2$}
  \label{fig:TwoD_Cost}
\end{subfigure}
\caption{Numerical results in two dimensions.}
\label{fig:TwoDResults}
\end{figure}

\begin{figure}
\centering
\begin{subfigure}{.5\textwidth}
  \centering
  \includegraphics[width=.8\linewidth]{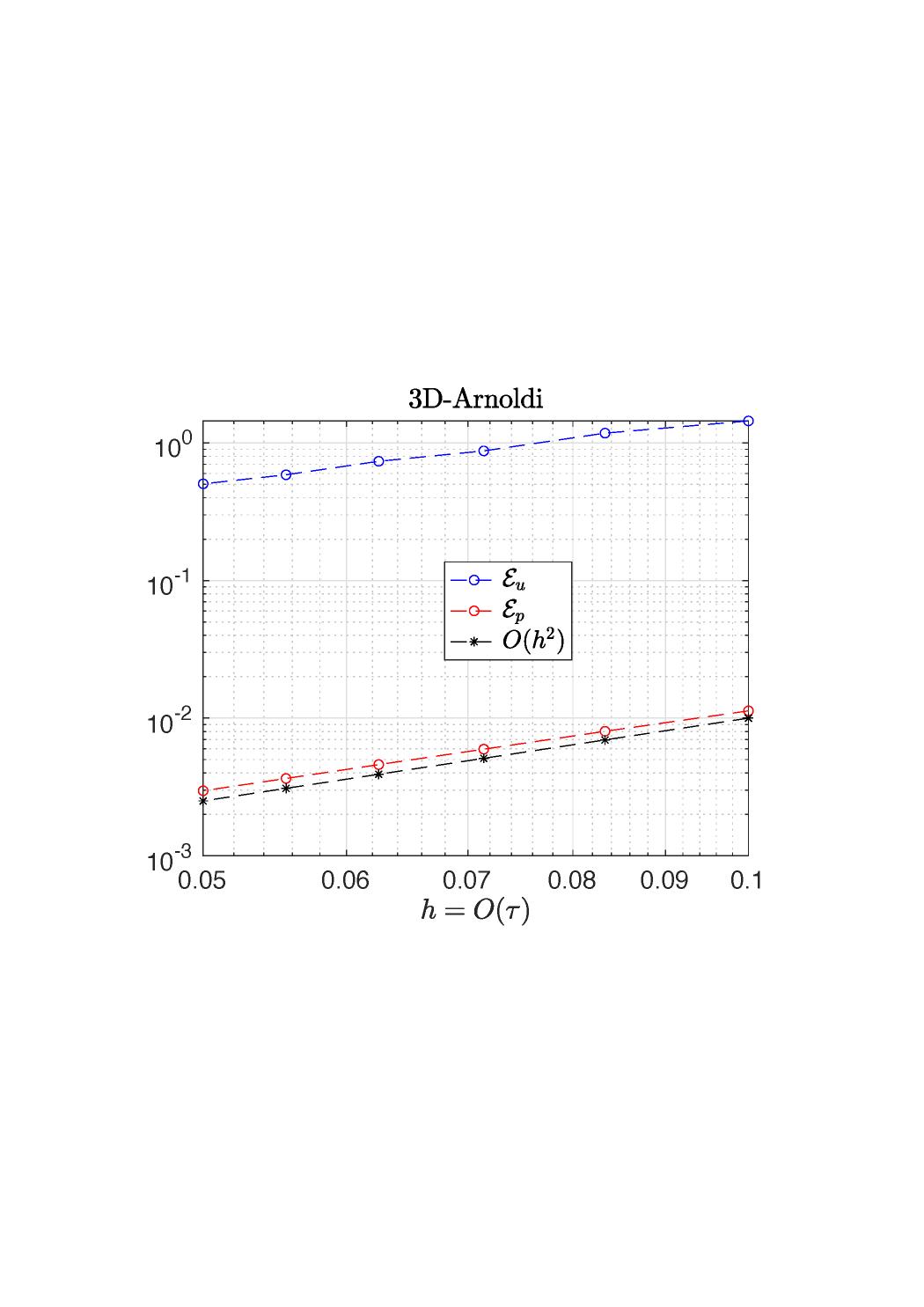}
  \caption{Convergence of $u$ and $p$ in $d=3$}
  \label{fig:ThreeD_Ex_1}
\end{subfigure}%
\begin{subfigure}{.5\textwidth}
  \centering
  \includegraphics[width=.8\linewidth]{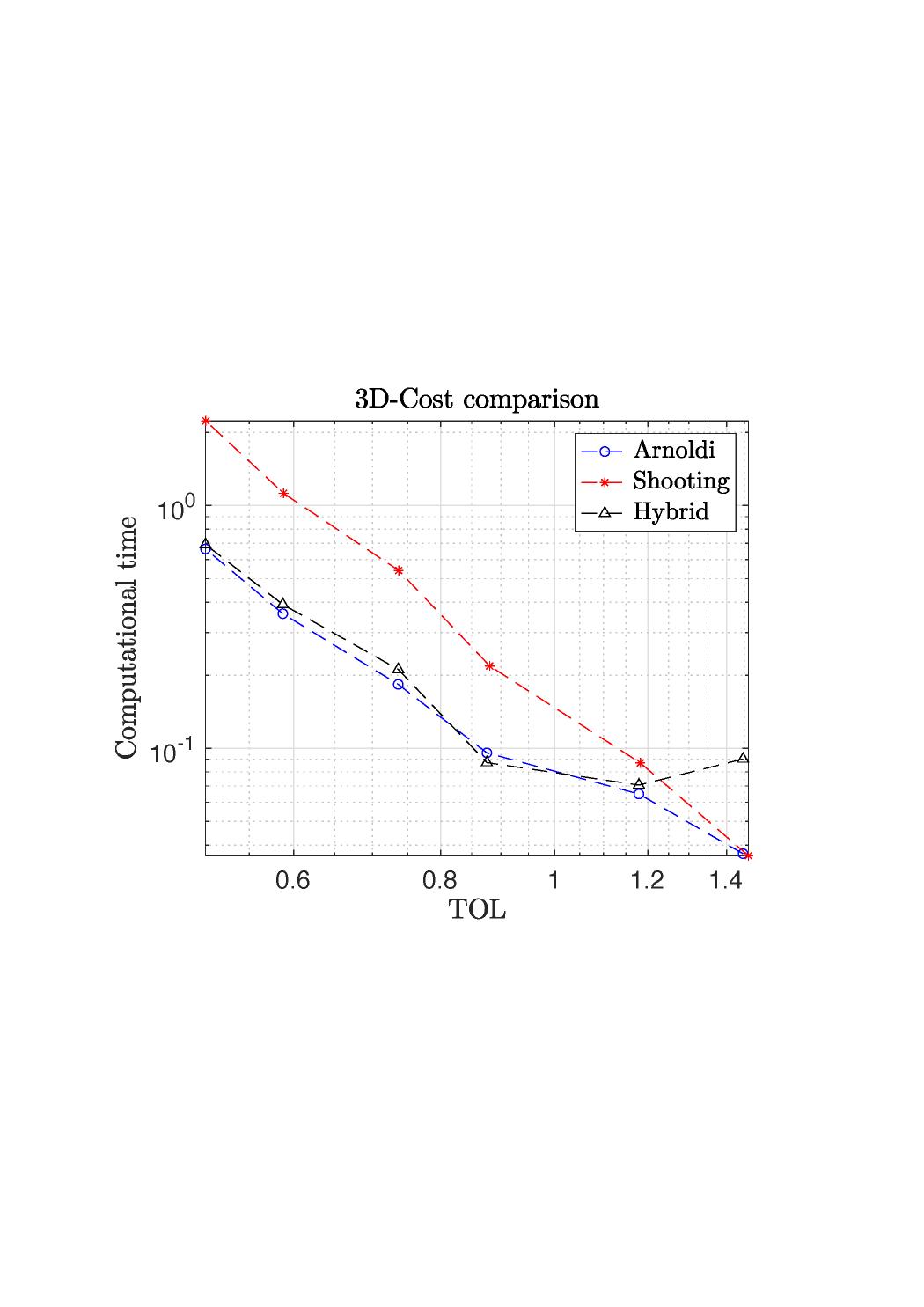}
  \caption{Cost comparison in $d=3$}
  \label{fig:ThreeD_Cost}
\end{subfigure}
\caption{Numerical results in three dimensions.}
\label{fig:ThreeDResults}
\end{figure}

In  this section, we consider the inverse problem of  determining the pair of solutions $\{ p,u\}$ in dimensions $d=1,2,3$: 
\begin{eqnarray*}
\partial_{t} u(t,x) & - & \Delta u(t,x) = f(t,x) + p(x), \quad \text{in }  K = (0,1)^d \times (0,T), \\
u(t,x) & = & 0, \quad \text{on }  \partial K \times (0,T),  \\
u(0,x) & = & 0, \quad \text{ in } K, \\
u(T,x) & = & \varphi(x), \quad \text{ in } K.
\end{eqnarray*}
To study the convergence properties, we assume an exact solution $u$ (only for the sake of comparison) of the form
\[
u(t,x) = \left( e^{-t} - 1 \right) \prod_{j=1}^{d} \sin^2(2 \pi x_j),
\]
and
\begin{eqnarray*} 
f(t,x) & = & - e^{-t} \left(\prod_{j=1}^{d} \sin^2(2 \pi x_j) + 8 \pi^2 \sum_{j=1}^{d} \cos(4 \pi x_j) \prod_{k\neq j}^{d} \sin^2(2 \pi x_k) \right), \\
p(x) & = & 8 \pi^2 \sum_{j=1}^{d} \cos(4 \pi x_j) \prod_{k\neq j}^{d} \sin^2(2 \pi x_k), \\
\varphi(x) & = & u(T,x). 
\end{eqnarray*}
In Figures \ref{fig:OneD_Ex_1},\ref{fig:TwoD_Ex_1},\ref{fig:ThreeD_Ex_1}, the relative errors 
\begin{align*}
    \mathcal{E}_u := \max_{k,j} \left| \dfrac{ u(t_k,x_j) - u_{k,j} }{u(t_k,x_j)} \right|, \quad  \mathcal{E}_p := \max_{j} \left| \dfrac{p(x_j) - p_{j}}{p(x_j)} \right|
\end{align*}
are depicted for dimensions $d=1,2,3$ respectively. Here $u_{k,j}$ and $p_j$ are numerical approximations to $u(t_k,x_j)$ and $p(x_j)$ respectively. The error plots include only the error corresponding the full Arnoldi approach. This is intentional since the error plots for the hybrid approach as well as the shooting algorithm are almost identical, and all show second order convergence rates in time and space. Note that in all of the simulations the time-step is simultaneously refined (proportional to the spatial stepsize). Moreover, the final time $T$ is set as $T=0.1$, and the number of basis vectors in the Arnoldi algorithm is chosen as $k=N$ in all of the simulations.  

To compare the efficiency of the methods, we also report, in Figures \ref{fig:OneD_Cost}, \ref{fig:TwoD_Cost},\ref{fig:ThreeD_Cost}, the actual computational time (measured in seconds) to reach a desired error tolerance. For dimensions $d=2,3$, the results clearly demonstrate the advantage of using the hybrid and the full Arnoldi approaches over the classical shooting method, whereas the shooting method seems to be more efficient only in the one-dimensional setting.  


\clearpage

\bibliographystyle{plain}
\bibliography{Ref}

\end{document}